\documentclass[11pt]{amsart}
\usepackage{amssymb}
\usepackage[hidelinks]{hyperref}

\newtheorem{theorem}{Theorem}[section]
\newtheorem{corollary}[theorem]{Corollary}
\newtheorem{proposition}[theorem]{Proposition}
\newcommand{\F}{\mathbb{F}_5}
\newcommand{\proofstep}[1]{\par\smallskip\noindent\textit{#1.}\enspace\ignorespaces}

\title[Normal abelian subgroups of finite \(5\)-groups]
      {Normal abelian subgroups of order \(5^6\) in finite \(5\)-groups}
\author{L\"{u} Gong}
\author{Jiaye Chen}
\author{Baojun Li}
\address{School of Mathematics and Statistics, Nantong University,
Nantong, Jiangsu 226019, China}
\email[L\"{u} Gong]{lieningzai1917@126.com}
\email[Jiaye Chen]{Cjyi@stmail.ntu.edu.cn}
\email[Baojun Li]{bjli@ntu.edu.cn}
\thanks{Correspondence should be addressed to Baojun Li at
\href{mailto:bjli@ntu.edu.cn}{\nolinkurl{bjli@ntu.edu.cn}}.}
\expandafter\providecommand\csname subjclassname@2020\endcsname
  {\textup{2020} Mathematics Subject Classification}
\date{}

\begin{document}

\begin{abstract}
A finite \(5\)-group containing an abelian subgroup of order \(5^6\)
has a normal abelian subgroup of that order. If the given subgroup is
elementary abelian, the normal subgroup can be chosen elementary
abelian. This resolves \(n=6\), the last open case of both parts of
Problem~21.51 of the \emph{Kourovka Notebook} for \(p=5\).
Together with earlier results, it follows that \(1\le n\le6\) is
the exact range in which each assertion holds for every finite
\(5\)-group at subgroup order \(5^n\).
In a minimal counterexample, five conjugate abelian subgroups generate
a group of class at most three. In class two we construct invariant
isotropic subgroups for the commutator pairing. In class three the
Jacobi identity produces a nonzero plane cubic over the field of five
elements. The five conjugate subgroups give five distinct singular
points spanning the plane, which forces the cubic to vanish.
\end{abstract}

\maketitle

\section{Introduction}

Let \(P\) be a finite \(p\)-group. If \(P\) contains an abelian
subgroup of order \(p^n\), for which values of \(n\) must it contain
a normal abelian subgroup of the same order? George Glauberman posed
this question and its elementary abelian counterpart as
Problem~21.51 of the \emph{Kourovka Notebook}
\cite[Problem~21.51, p.~169]{KN}. The second question asks whether,
when the given subgroup is elementary abelian, the normal subgroup
can also be chosen elementary abelian.

For odd \(p\), Konvisser and Jonah proved both assertions for
\(n\le 5\) using a projective approach
\cite[Theorems~4.1 and~5.5]{KJ}. Alperin and Glauberman developed
methods for finding invariant abelian subgroups, including a
construction of invariant elementary abelian subgroups from suitable
families of conjugates \cite[\S\S2, 4]{AG}. Glauberman subsequently
proved that both assertions hold for \(n<(p+2)/2\)
\cite{GlExist}, giving a range that grows with \(p\).

In the other direction, Berkovich constructed a group of order
\(5^{11}\) containing an elementary abelian subgroup of order \(5^7\)
but no normal abelian subgroup of that order
\cite[Theorem~14.4]{Ber}. Glauberman described this example and
constructed counterexamples for primes \(p\ge 7\) at
\(n=(p+9)/2\) \cite{Gl}. Direct products with elementary abelian
groups extend these counterexamples to larger \(n\).
Thus, for \(p=5\), the preceding results establish both assertions
through \(n=5\) and rule out both for \(n\ge 7\). The following theorem
settles the remaining case \(n=6\) in both parts of Glauberman's
problem.

\begin{theorem}\label{thm:six}
Let \(P\) be a finite \(5\)-group.
\begin{enumerate}
\item[\textup{(a)}] If \(P\) has an abelian subgroup of order \(5^6\),
then \(P\) has a normal abelian subgroup of order \(5^6\).
\item[\textup{(b)}] If \(P\) has an elementary abelian subgroup of
order \(5^6\), then \(P\) has a normal elementary abelian subgroup of
order \(5^6\).
\end{enumerate}
\end{theorem}

\begin{corollary}\label{cor:threshold}
For each of the two assertions in Theorem~\ref{thm:six}, the positive
integers \(n\) for which the assertion with \(5^n\) in place of
\(5^6\) holds for every finite \(5\)-group are precisely
\(1,2,3,4,5,6\).
\end{corollary}

Why does a plane cubic enter a question about commuting elements?
In a smallest counterexample one can find an abelian subgroup of
order \(5^6\) with exactly five conjugates. In the last possible case,
their pairwise commutators generate a three-dimensional vector space
\(S\). Each subgroup determines a plane in \(S\), hence a point of
the dual projective plane. Commutators of weight three satisfy the
Jacobi identity, which makes the coefficients of a trilinear form
symmetric. Its diagonal is a nonzero cubic. Commutativity within
each conjugate says that all first derivatives vanish at the
corresponding point. The five points are distinct and span the plane.
A line through two of them already accounts for four zeros counted
with multiplicity on a cubic, so that line must be a component.
There are at least four such lines.

\begin{proof}[Deduction of Corollary~\ref{cor:threshold}]
For \(1\le n\le5\) the elementary abelian and abelian cases follow
from Konvisser and Jonah \cite[Theorems 4.1 and 5.5]{KJ},
respectively. Theorem~\ref{thm:six} supplies \(n=6\).
Berkovich constructed a group \(G\) of order \(5^{11}\) with an
elementary abelian subgroup of order \(5^7\) and no normal abelian
subgroup of that order \cite[Theorem 14.4]{Ber}. Glauberman describes this construction in
\cite[pp.~61--62]{Gl}.
For \(n\ge7\) take \(G\times C_5^{\,n-7}\), where \(C_5\) is cyclic
of order five. The indicated elementary abelian subgroup of \(G\)
and the second factor generate an elementary abelian subgroup of
order \(5^n\). A normal abelian subgroup of order \(5^n\) would
project to a normal abelian subgroup of \(G\) of order at least
\(5^7\). Every normal subgroup of a finite \(5\)-group admits a
series of subgroups normal in the whole group with factors of order
five. The image would therefore contain a normal abelian subgroup of
\(G\) of order \(5^7\), a contradiction.
\end{proof}

Section~\ref{sec:cubic} gives this argument assuming two group-theoretic
reductions. Section~\ref{sec:orbit} proves the first by bounding the
class of the group generated by the five conjugates.
Sections~\ref{sec:relative-class-two} and~\ref{sec:ordinary-class-two}
prove the second, a replacement theorem for groups of class at most
two. Its elementary abelian case uses polynomial interpolation over
\(\F\) and an extremal choice of an isotropic subspace, in the spirit
of the method of Alperin and Glauberman \cite[\S\S2, 4]{AG}.

\subsection*{Computational exploration}
Exploratory JavaScript calculations over \(\F\) used row reduction
and subspace enumeration to test commutator and Jordan block
configurations. The arguments below establish the required identities
and dimension bounds directly.

\subsection*{AI Declaration}
We used OpenAI Codex (GPT-6.2) to revise the exposition, check steps in
the proofs and bibliographic references, and run computations on finite
groups. We also used Codex to prepare the journal version. We take
full responsibility for the manuscript, including its mathematical
arguments.

\section{Five conjugates and a plane cubic}\label{sec:cubic}

Fix one assertion of Theorem~\ref{thm:six}. In case~(a) call a
subgroup \emph{admissible} if it is abelian. In case~(b) require it
to be elementary abelian. Suppose that \(P\) is a counterexample of
least order. Write \(\lambda(A)=\log_5|A|\) for a finite \(5\)-group
\(A\). Every admissible normal subgroup of \(P\) has length at most
five. Indeed, a nontrivial normal subgroup meets the centre of a finite
\(5\)-group nontrivially. Induction in successive quotients gives a
series normal in \(P\) with factors of order five, and admissibility
passes to subgroups. All commutators are left-normed. All vector
spaces in the proof are over \(\F\).

\begin{proposition}\label{prop:orbit}
There are five distinct conjugate admissible subgroups
\(B_i\) of order \(5^6\), indexed by \(i\in\mathbb Z/5\mathbb Z\),
and an element \(g\in P\) with \(B_i^g=B_{i+1}\).
Their product \(Q\) is normal in \(P\), each \(B_i\) is normal in
\(Q\), and \(Q\) has class at most three. In case~\textup{(b)},
\(Q\) has exponent five.
\end{proposition}

\begin{proposition}\label{prop:relative}
Let a finite \(5\)-group \(H\) act by automorphisms on a finite
\(5\)-group \(T\) of class at most two. If \(T\) contains an
admissible subgroup of order \(5^6\), it contains an \(H\)-invariant
admissible subgroup of that order, for either choice of admissibility.
\end{proposition}

\begin{proof}[Proof of Theorem~\ref{thm:six}]
Apply Proposition~\ref{prop:orbit}. If a subgroup \(T\le Q\) normal
in \(P\) has class at most two and contains an admissible subgroup
of order at least \(5^6\), Proposition~\ref{prop:relative}, applied
to the conjugation action of \(P\), contradicts our choice of \(P\).
In particular \(Q\) has class exactly three. In case~(b) every
abelian subgroup of \(Q\) is elementary abelian by the exponent
assertion of Proposition~\ref{prop:orbit}.

\proofstep{The centre and the derived group}
Put \(Z=Z(Q)\), \(D=Q'Z\), and \(C=[D,Q]=\gamma_3(Q)\). Then \(D\) is abelian, \(C\le Z\),
and all three groups are characteristic in \(Q\). Thus \(\lambda(D)\le5\). We claim that
\(C\) is noncyclic. Indeed, if \(C\) is cyclic, take \(E=C_Q(D)\). This is characteristic
in \(Q\) and has class at most two, since \(E'\le D\le Z(E)\). The central commutator map
gives an injection
\[
 B_0/(B_0\cap E)\ \hookrightarrow\
 \operatorname{Hom}\bigl(D/(B_0\cap D),C\bigr).
\]
For a cyclic target the number of homomorphisms is at most the order of the domain, by its
cyclic decomposition. Thus \(|B_0:B_0\cap E|\le |D:B_0\cap D|\). It follows that the
abelian subgroup \((B_0\cap E)D\) of \(E\) has order at least \(|B_0|=5^6\), contrary to
Proposition~\ref{prop:relative}. In particular \(Z\) contains a noncyclic subgroup and
\(\lambda(Z)\ge2\).

\proofstep{The middle layer}
Put \(S=D/Z\), \(V=Q/D\), and \(S_i=((B_i\cap D)Z)/Z\). At this stage \(S,V\) are finite
abelian groups. The commutator modulo \(Z\) gives an alternating biadditive map with
values generating \(S\). It is well defined since \([D,Q]\le Z\). We have
\[
 S=\sum_{i<j}\bigl([B_i,B_j]Z/Z\bigr),
 \qquad [B_i,B_j]Z/Z\ \le S_i\cap S_j .
\]
Here we use that \(Q\) is generated by the five \(B_i\), each abelian and normal in \(Q\).
If some \(S_i=S\), conjugacy gives \(S_j=S\) for every \(j\). Then for each \(j\), the
equality \(D=(B_j\cap D)Z\) shows that \(D\) centralises \(B_j\), contrary to class three.
Thus each \(S_i\) is proper in \(S\). Also \(S\ne0\). If the common length of the \(S_i\)
is at most one, their orbit under \(g\) is either one subgroup or five distinct subgroups
of order five. In the first case their pair intersections lie in a fixed proper subgroup.
In the second those intersections vanish. Both contradict the displayed generation
property. Hence \(\lambda(S_i)\ge2\) and \(\lambda(S)\ge3\). Since
\(\lambda(D)=\lambda(S)+\lambda(Z)\le5\), we conclude
\[
 \lambda(Z)=2,\qquad \lambda(S)=3,\qquad \lambda(D)=5,\qquad \lambda(S_i)=2 .
\]
Moreover \(C=Z\) is elementary abelian of order \(25\), because \(C\) is noncyclic. Since
\([D,Q]\le Z\) is central of exponent five, \([d^5,q]=[d,q]^5=1\). Hence \(D^5\le Z\) and
\(S\) is an elementary abelian group of rank three.

\proofstep{The evaluation pairing}
We next show that \(C_Q(D)=D\). Otherwise \(D\) together with an element outside \(D\) in
its centraliser generates an abelian group of order at least \(5^6\) inside the
characteristic subgroup \(C_Q(D)\) of class at most two. This contradicts
Proposition~\ref{prop:relative}. The central commutator gives an injective homomorphism
\[
 \theta:V\lhook\joinrel\longrightarrow\operatorname{Hom}(S,Z),\qquad
 \theta(v)(s)=[v,s],
\]
using lifts from the indicated quotients. The defining bracket is independent of the lifts
because \(D\) is abelian and \(Z\) is central. Its kernel on \(Q\) is \(C_Q(D)\). Thus
\(V\) is elementary abelian as well, and the maps \(\theta\) and \(\alpha:\bigwedge^2 V\to
S\) are linear, with \(\alpha\) onto. For each \(i\) write \(S_i=\ker\ell_i\), where
\(\ell_i\in S^*\setminus0\), and \(V_i=B_iD/D\). Since \(B_i\) centralises \(B_i\cap D\),
we have \(\theta(V_i)\subseteq Z\otimes\langle\ell_i\rangle\), a space of dimension two.
Meanwhile \(\lambda(B_i\cap D)\le\lambda(Z)+\lambda(S_i)=4\), so \(\dim V_i\ge2\) and
equality holds throughout. The orbit of \(S_0\) has size one or five. In the former case
every pair intersection lies in the same proper plane, contrary to the generation of
\(S\). Thus the planes are all distinct. Their defining functionals span \(S^*\). Indeed,
if their span had dimension two, every pair of planes would meet in the same line,
contradicting the generation of \(S\) by pair intersections. Thus the \(V_i\) together
span all of \(Z\otimes S^*\), and
\[
 V=\operatorname{Hom}(S,Z),\qquad
 V_i=Z\otimes\langle\ell_i\rangle,\qquad |Q|=5^{11}.
\]

\proofstep{Jacobi gives a cubic}
For \(x,y,z\in V\) the class-three commutator identity is
\[
 \theta(x)\alpha(y,z)+\theta(y)\alpha(z,x)+\theta(z)\alpha(x,y)=0.
\]
This is the weight-three Hall--Witt identity. Equivalently apply the Jacobi identity in
the Lie ring of the lower central series. The added central terms in forming \(V\) and
\(S\) have zero bracket, so the identity descends for either choice of admissibility.
Choose a basis \(e_1,e_2\) of \(Z\) and write \(X_a=e_1\otimes a\), \(Y_b=e_2\otimes b\)
for \(a,b\in S^*\). The \(e_2\) component of Jacobi on \(X_a,X_c,Y_b\) says
\(b(\alpha(X_a,X_c))=0\) for all \(b\). Thus \(\alpha(X_a,X_c)=0\), and exchanging the two
basis vectors gives \(\alpha(Y_a,Y_c)=0\).
Define
\[
 F(a,b,c)=c\bigl(\alpha(X_a,Y_b)\bigr).
\]
The other component of the same identity gives \(F(c,b,a)=F(a,b,c)\). Taking the \(e_2\)
component of Jacobi on \(X_a,Y_b,Y_c\) gives \(F(a,c,b)=F(a,b,c)\). Thus \(F\) is
symmetric. It is nonzero since \(\alpha\) is surjective and vanishes on both unmixed
pairs. As \(6\) is invertible in \(\mathbb F_5\), polarisation gives a nonzero homogeneous
cubic \(f(a)=F(a,a,a)\).

Each \(V_i\) is isotropic for \(\alpha\), because its lifts in \(B_i\) commute.
Consequently \(F(\ell_i,\ell_i,c)=0\) for every \(c\in S^*\). In particular
\(f(\ell_i)=0\), and
\[
  (df)_{\ell_i}(c)=3F(\ell_i,\ell_i,c)=0 .
\]
Thus each \([\ell_i]\) is a singular point of the cubic \(f=0\) in \(\mathbb{P}(S^*)\).
These five points span the plane. For any two of these points, the restriction of \(f\) to
their line vanishes identically, since a binary cubic with two distinct double roots is
zero. Thus every secant of the five points divides \(f\) as a linear factor. Choose three
noncollinear points among the five. Their sides give three distinct secants. A fourth
point lies on at most one of these sides, so joining it to a suitable vertex gives a
fourth secant. A nonzero cubic cannot have four distinct linear factors. This
contradiction proves Theorem~\ref{thm:six}.
\end{proof}

\section{Bounding the class of the normal closure}\label{sec:orbit}

\begin{proof}[Proof of Proposition~\ref{prop:orbit}]
Take a maximal subgroup \(M<P\) containing an admissible witness of
order \(5^6\). The witness is proper since \(P\) is a counterexample.
Minimality applied to \(M\) gives an admissible subgroup
\(B\trianglelefteq M\) of order \(5^6\). Its normaliser in \(P\)
equals \(M\). Fix \(g\in P\setminus M\) and let
\[
 B_i=B^{g^i}\quad (i\in\mathbb Z/5\mathbb Z),\qquad
 Q=\langle B_0,B_1,B_2,B_3,B_4\rangle .
\]
Here \(M\trianglelefteq P\), \(g^5\in M\), and the \(B_i\) are
distinct, normal in \(M\), and contained in \(M\). In particular they
are normal in \(Q\), and \(Q\trianglelefteq P\).

\proofstep{Intersection subgroups}
For a nonempty set of labels \(I\), put \(B_I=\bigcap_{i\in I}B_i\).
Let \(L_r\) be the product of all \(B_I\) with \(|I|=r\), and put
\(L_r=1\) for \(r>5\). The factors are normal in \(M\), so these
products are subgroups. They are normal in \(P\) because conjugation
permutes the labels. The subgroups form a descending series with
\(L_1=Q\). If \(I\cap J\ne\varnothing\), then \(B_I,B_J\) lie in a
common abelian group \(B_i\) and commute. If \(I\cap J=\varnothing\),
their normality in \(M\) gives
\([B_I,B_J]\le B_{I\cup J}\). Consequently
\[
 [L_r,L_s]\le L_{r+s},\qquad \gamma_r(Q)\le L_r .
\]
In particular \(L_3\) is admissible. In part (b) this follows because
its elementary abelian factors commute. Thus \(|L_3|\le 5^5\).

\proofstep{Fourfold intersections}
If \(L_4>L_5\), take \(x\) in a fourfold intersection and outside the
remaining \(B_j\). The elements \(x^{g^s}\), \(0\le s<5\), have
different missing labels. Under the homomorphism
\[
 L_4\longrightarrow\prod_i L_4/(L_4\cap B_i)
\]
each maps to an element supported at its own missing label, of order
at least five. The kernel is \(L_5\). Hence \(|L_4/L_5|\ge5^5\).
As \(L_4\le L_3\), this forces \(L_3=L_4\) and \(L_5=1\).
It follows that \([L_3,Q]=1\), so \(\gamma_4(Q)=1\).

We may therefore assume \(L_4=L_5\). This subgroup is central in
\(Q\), since it lies in every \(B_i\). Thus \(\gamma_5(Q)=1\).
In both cases considered so far the class is less than five. In
part (b), the regularity theorem for finite \(p\)-groups of class
less than \(p\) implies that the elements of order dividing five
in \(Q\) form a subgroup (see \cite[\S4]{AG}). They include all the \(B_i\),
so \(Q\) has exponent five.

\proofstep{Exclusion of class four}
Suppose \(\gamma_4(Q)\ne1\). Since \(\gamma_5(Q)=1\), the fourth
lower central factor is generated by left-normed commutators of four
elements from the union of the \(B_i\). This follows by collecting
commutators modulo the next lower central term at each step.
Choose one such nontrivial commutator
\[
 [t,b]\ne1,\qquad t=[a_1,a_2,a_3].
\]
The four labels must be distinct. A left-normed commutator lies in
every factor whose label has occurred, so a repeated label makes it
trivial. Let \(I\) be the triple of labels for \(t\), and \(j\) the
label for \(b\). For \(0\le s<5\), write \(t_s=t^{g^s}\),
\(b_s=b^{g^s}\) and define
\[
 f_s:L_3/L_4\longrightarrow L_4,\qquad xL_4\longmapsto[x,b_s].
\]
The maps are homomorphisms because \([L_3,Q]\le L_4\le Z(Q)\).
We have \(f_s(t_sL_4)\ne1\). Also \(f_s(t_rL_4)=1\) if
\(j+s\in I+r\). For fixed \(s\), the three elements of \(I\) exclude
three possible columns \(r\). The only possible nontrivial columns
are \(s\) and \(s+\delta\), with the same nonzero
\(\delta\in\mathbb Z/5\mathbb Z\) for every row.

In the order \(0,\delta,2\delta,3\delta\), each new column
\(t_{k\delta}L_4\) is outside the subgroup generated by its
predecessors. Indeed, \(f_{k\delta}\) is nontrivial on this column and
trivial on all preceding ones. Each adjunction increases the order
by a factor of at least five. Thus \(|L_3/L_4|\ge5^4\). As
\([t,b]\ne1\), the bound \(|L_3|\le5^5\) now gives
\[
 |L_3|=5^5,\qquad |L_4|=5.
\]

We have \(C_Q(L_3)=L_3\). Otherwise \(C_Q(L_3)/L_3\) is a
nontrivial normal subgroup of \(P/L_3\), hence meets its centre in
a subgroup of order five. The inverse image of this subgroup is
normal in \(P\), has order \(5^6\), and is abelian because it is
generated by \(L_3\) and one element centralising \(L_3\). In
part (b) it is elementary abelian by the exponent calculation above.
This contradicts our choice of \(P\).

Finally put \(U=B\cap L_3\). The central commutator pairing induces an injection
\[
 B/U\hookrightarrow \operatorname{Hom}(L_3/U,L_4),\qquad
 bU\longmapsto \bigl(xU\mapsto[x,b]\bigr).
\]
The pairing vanishes on \(U\) in each argument since \(B\) and \(L_3\)
are abelian. Its kernel on \(B\) is \(U\) since \(C_Q(L_3)=L_3\).
As \(L_4\) has order five, the target has order at most \(|L_3/U|\). This would give
\(|B|\le|L_3|=5^5\), whereas \(|B|=5^6\). This excludes class four and proves the proposition.
\end{proof}

\section{Invariant elementary abelian subgroups in class two}
\label{sec:relative-class-two}

\begin{proof}[Proof of Proposition~\ref{prop:relative} in the elementary abelian case]
We work over \(k=\mathbb F_5\). In a group of class at most two,
the elements of order dividing \(5\) form a characteristic subgroup
of exponent dividing \(5\). Indeed, if \(x^5=y^5=1\), then
\([x,y]^5=[x^5,y]=1\), and the class-two power formula gives
\((xy)^5=1\), since \(5\mid\binom{5}{2}\).
Replacing \(T\) by this characteristic subgroup preserves the witness
and the action. We can therefore assume \(T\) has exponent \(5\).
Write \(W=T'\), \(V=T/W\), \(m=\dim W\), and \(d=\dim V\).
The commutator induces an \(H\)-equivariant surjection
\[
  \beta:\bigwedge\nolimits^2 V\longrightarrow W.
\]
Since \(W\) is central, the product of the given elementary abelian
subgroup with \(W\) is elementary abelian. If \(m<6\), it follows that
there is a totally isotropic subspace of dimension \(r=6-m\) for
\(\beta\). An \(H\)-invariant totally isotropic \(r\)-space lifts to
an \(H\)-invariant elementary abelian subgroup of order \(5^6\).
Throughout, a finite \(5\)-group acting linearly over \(k\) fixes a
nonzero vector on each nonzero module and preserves a complete flag.
This follows by counting orbits and then passing to successive quotients.

If \(m\le 2\), Theorem~2.3 of Konvisser and Jonah \cite{KJ} shows that
the number of elementary abelian subgroups of \(T\) of order \(5^6\)
is \(1\) modulo \(5\). The action of \(H\) on this set has a fixed
point. If \(m\ge 6\), use a six-dimensional member of an invariant
flag in \(W\). If \(m=5\), use an invariant line in \(V\) and take its
inverse image. It remains to treat \(m=3,4\).

\proofstep{An interpolation criterion}
For \(h\in H\) let \(h_V,h_W\) denote its actions on the two spaces.
Suppose that \(h_V^5=h_W^5=1\), and for \(M=V,W\) put
\[
  h_M(t)=\sum_{i=0}^4\binom ti(h_M-1)^i .
\]
The polynomials are defined over \(k\), and
\(h_M(t+u)=h_M(t)h_M(u)\) by the binomial convolution identity and
\((h_M-1)^5=0\). At every \(i\in k\), these matrices are the powers
\(h_M^i\). Consequently if
\(\beta\circ\bigwedge^2 h_V(t)\) has degree at most four, then
\[
    \beta\circ\bigwedge\nolimits^2 h_V(t)=h_W(t)\circ\beta
\]
identically, by comparison at the five elements of \(k\).

If every \(h\in H\) satisfies these hypotheses, an \(H\)-invariant
isotropic \(r\)-space exists whenever an ordinary one does.
To see this, let \(X\) be the closed variety of isotropic \(r\)-spaces
in the Grassmannian. The one-parameter subgroups \(h_V(t)\) preserve \(X\).
In a basis adapted to an \(H\)-invariant complete flag, they are
upper unitriangular. The induced matrices on the exterior power
are also upper unitriangular in the lexicographically ordered wedge
basis. Among the \(k\)-rational points of \(X\), choose one with
least possible index of the last nonzero Pl\"ucker coordinate, and
normalise that coordinate to one. If some \(h\) moves the point,
the coordinates of its image under \(h_V(t)\) are polynomials, at
least one of them nonconstant. The distinguished coordinate remains
one and all later coordinates remain zero. The vector of coefficients
of highest degree is nonzero and supported strictly earlier.
Its projective point still belongs to \(X(k)\). Indeed, substitution into each homogeneous
equation of \(X\) gives a polynomial identity, whose leading coefficient gives the same
equation for this vector.
This contradicts minimality. The same chosen point is therefore
fixed by every element of \(H\).

On two distinct Jordan blocks of lengths \(a,b\), the degree of
\(\bigwedge^2 h_V(t)\) is at most \(a+b-2\).
Within a single block of length \(a\ge2\), it is at most \(2a-4\).
Indeed, writing the basis as \(x_i=(h_V-1)^ix_0\), each term mapping
\(x_i\wedge x_j\) to \(x_s\wedge x_t\) has degree at most
\(s+t-i-j\), with \(i\ne j\) and \(s\ne t\).

\proofstep{Three-dimensional commutator space}
Here \(r=3\). First observe that it is enough to find an \(H\)-fixed
nonzero \(v\in V\) such that
\[
    d-1-\operatorname{rank}\beta(v,-)\ge4.
    \tag{*}
\]
Indeed \(C=\ker\beta(v,-)\) is invariant and the form induced on
\(C/kv\) takes values in \(W\). Choose an invariant four-space \(L\)
in this quotient. The kernel \(K\) of \(\bigwedge^2 L\to W\) has
dimension at least three. After choosing a volume form on \(L\),
\(z\mapsto z\wedge z\) gives an \(H\)-invariant scalar quadratic
form on \(K\), since \(H\) has determinant one on \(L\).
This form has a nonzero isotropic vector, since every quadratic form of dimension at
least three over \(k\) is isotropic. For instance, in the nondegenerate
diagonal ternary case \(ax^2+by^2+cz^2\), the sets of values
\(ax^2\) and \(-c-by^2\) have three elements each and intersect.
There is also a fixed nonzero zero. If the radical is nonzero take a
fixed vector in it. Otherwise, if the fixed subspace were anisotropic,
it would split off orthogonally. Any nonzero orthogonal complement
would itself have a nonzero fixed vector. This forces the entire
space to be fixed, in which case any isotropic vector suffices.
The nonzero zeros in \(\bigwedge^2 L\) are precisely the decomposable
bivectors, by the rank criterion for alternating matrices of size
four. A fixed such vector yields an invariant isotropic two-plane
in \(L\), whose inverse image under \(C\to C/kv\) is the desired three-plane.
In particular, \((*)\) always applies for \(d\ge8\).

Suppose \(d\le7\). Fix an \(H\)-invariant complete flag
\(0=F_0<F_1<\cdots<F_d=V\). If some \(h_V\) has a block of length
at least six, write \(D=h_V-1\). Since \(h_W^5=1\) and
\(h_V^5=1+D^5\), equivariance gives
\(\beta(D^5 V,\ker D^5)=0\). Every element of \(D^5 V\) thus
has contraction rank at most \(\operatorname{rank}D^5\).
For a single block of length \(d=6\) or \(7\), \(F_1=D^{d-1}V\)
and this bound is \(d-5\), which gives \((*)\). The only remaining
type with a long block is \((6,1)\). Since \(h\) acts trivially on the factors of the
flag, \(DF_i\subseteq F_{i-1}\). Hence \(D^5V\) is a line in \(F_2\).
If this line is \(H\)-invariant, it provides a fixed vector of contraction rank at most
one. Otherwise two distinct \(H\)-translates span \(F_2\).
Contraction ranks are subadditive and each translate consists of
vectors of rank at most one. A nonzero vector in \(F_1\) now has
rank at most two, giving \((*)\).

It remains to verify the interpolation criterion when all block lengths
are at most five. Write \(D=h_V-1\), \(E=h_W-1\), so that \(E^3=0\).
For two block generators \(x,y\) of lengths \(a,b\), respectively,
\[
 \beta(D^i x,D^j y)=0\qquad\text{if }i+j\ge\min(a,b)+2.
\]
Indeed,
\(\beta(Dx,y)=E\beta(x,h_V^{-1}y)
+\beta(x,(h_V^{-1}-1)y)\).
The two operations on the right commute. After \(i\) iterations every
surviving term contains at most two factors \(E\), and therefore at
least \(i-2\) additional factors \(D\) on the second variable, since
\(h_V^{-1}-1=-h_V^{-1}D\). This proves the bound using \(b\). Interchanging the variables gives the bound using \(a\).
For \(d\le7\), the only distinct-block pairs with an uncorrected
degree above four have lengths \((5,2)\) and \((4,3)\). All contributions
of degree above four in those pairs vanish by this bound.

Within a block of length five, a term of degree at least five has
distinct output indices of sum at least six. Thus it suffices to
prove \(\beta(x_2,x_4)=\beta(x_3,x_4)=0\). In the cyclic basis
\(v_i=h_V^i x_0\), \(i\in k\), the coordinates of \(x_4\) are all
one. Those of \(x_2,x_3\) are, respectively,
\((i+1)(i+2)/2\) and \(-i-1\), as follows by expanding \(D^j x_0\).
For either coordinate polynomial \(f\) and for each shift \(s=1,2\),
the coefficient of \(\beta(v_i,v_{i+s})\) in the bracket is
\(f(i)-f(i+s)\), of degree at most one. These two shifts account
for all unordered pairs. Since
\(\beta(v_i,v_{i+s})=(1+E)^i\beta(v_0,v_s)\), each weighted sum
of coefficients is a sum over \(i\in k\) of a polynomial of degree at
most \(1+2=3\), hence vanishes in \(k\). Within blocks of smaller
length the uncorrected exterior degree is already at most four.
The interpolation criterion proves the claim for \(m=3\).

\proofstep{Four-dimensional commutator space}
Here \(r=2\). If \(d\ge6\), take a nonzero fixed vector \(v\). Then
\(\ker\beta(v,-)\) has dimension at least two, and a fixed line in
\(\ker\beta(v,-)/kv\) gives the required invariant plane.
Assume \(d\le5\). If some \(h_V\) has a block of length five, it
has a single block, and the invariant complete flag of \(H\) is its
unique invariant complete flag. Its two-dimensional member is
\(\langle x_3,x_4\rangle\). In the cyclic-coordinate calculation
above the differences for \(x_3,x_4\) have degree zero, whereas
\(E^4=0\). The weighted coefficient sums still have degree at most
three and vanish. This invariant plane is therefore isotropic.
Finally, if all block lengths are at most four, a block of length
four has at most one other block, of length one. The distinct-block
degree bounds and the single-block bound are consequently all at
most four. Interpolation gives the required invariant plane.
\end{proof}

\section{Invariant abelian subgroups in class two}
\label{sec:ordinary-class-two}

\begin{proof}[Proof of Proposition~\ref{prop:relative} in the abelian case]
Put \(Z=Z(T)\), \(V=T/Z\), and \(w=\log_5|Z|\). We write the
abelian groups \(V\) and \(Z\) additively. The commutator induces an
\(H\)-equivariant alternating biadditive pairing
\[
 \beta:V\times V\longrightarrow Z.
\]
Its radical is zero by the definition of \(Z\). Hence
\(\exp V\leq\exp Z\), since multiplication by \(\exp Z\) sends
\(V\) into the radical. The image of an abelian witness in \(V\) is
totally isotropic and has order at least \(5^{6-w}\).

Every finite abelian \(5\)-group with an \(H\)-action has an
\(H\)-invariant series with successive factors of order five.
Indeed, its subgroup of elements killed by five has a nonzero fixed
vector if the group is nontrivial, and one can continue in the
quotient. In particular, an invariant abelian subgroup of order at
least \(5^6\) suffices. If \(w\geq6\), use \(Z\) itself. Otherwise
put \(k=6-w\). It suffices to find an invariant totally isotropic
subgroup of \(V\) of order at least \(5^k\), since its inverse image
in \(T\) is abelian. The case \(k=1\) follows from an invariant
subgroup of order five in \(V\).

\proofstep{Cyclic centre}
Suppose \(Z\) is cyclic. The pairing is perfect, since its adjoint
\(V\to\operatorname{Hom}(V,Z)\) is injective and both sides have
the same order by \(\exp V\leq\exp Z\). For every subgroup
\(U\leq V\) we have
\[
 |U^\perp|=|\operatorname{Hom}(V/U,Z)|=|V|/|U|.
\]
Indeed, the perfect adjoint identifies \(U^\perp\) with the
homomorphisms vanishing on \(U\). In particular every totally
isotropic \(U\) has order at most \(\sqrt{|V|}\).
We construct an invariant totally isotropic subgroup attaining this
bound. If \(V\neq0\), take an invariant subgroup \(L\) of order
five. It is cyclic and hence isotropic. The pairing on
\(L^\perp/L\) is nondegenerate because
\((L^\perp)^\perp=L\), as follows from the order formula. Its domain
has order \(|V|/25\). Induction, starting with the zero group,
gives an invariant isotropic subgroup there of order
\(\sqrt{|V|}/5\). Its inverse image in \(L^\perp\) attains the
required bound. Since there is already an isotropic witness of order
at least \(5^k\), this proves the claim for cyclic \(Z\).

\proofstep{Elementary abelian central quotient}
Suppose \(5V=0\). Then \(\beta\) takes values in
\(Z[5]=\{z\in Z:5z=0\}\). Enlarge this vector space to an
\(H\)-module \(W\) of dimension \(w\) by adjoining a trivial
summand. View \(\beta\) as a map into \(W\). On \(V\oplus W\)
define the multiplication
\[
 (v,z)(v',z')=(v+v',z+z'+\tfrac12\beta(v,v')).
\]
Bilinearity gives associativity. The resulting group has exponent
five and class at most two, and the diagonal action of \(H\) is by
automorphisms. Over the given isotropic witness in \(V\), the full
inverse image contains an elementary abelian subgroup of order
\(5^6\). The elementary abelian case of Proposition~\ref{prop:relative}, proved in
Section~\ref{sec:relative-class-two}, gives an \(H\)-invariant one
of this order. Its projection to \(V\) is invariant and isotropic,
with dimension at least \(6-\dim W=k\). This proves the claim in
this case.

\proofstep{Central quotients of higher exponent}
We may now suppose \(\exp V=5^e\) with \(e\geq2\), \(Z\) noncyclic,
and \(w\leq4\). Since \(\exp Z\geq25\), we have \(w=3\) or \(4\).
The highest-order cyclic direct factors of \(V\) have multiplicity
at least two. Otherwise, writing \(x\) for a generator of the unique
such factor, \(5^{e-1}x\) pairs trivially both with \(x\) and with
every lower-order factor. This would give a nonzero radical element.
It follows that \(R=5^{e-1}V\) is a characteristic elementary
abelian subgroup of rank at least two. Moreover
\(\beta(R,V[5])=0\) by biadditivity. Thus \(R\) is isotropic.
For \(w=4\), this gives the required order \(5^2\).

For \(w=3\), necessarily \(Z\cong C_{25}\times C_5\) and
\(e=2\). We need isotropic order \(5^3\). If \(R\) has rank at
least three we are done. If \(R\) has rank two and \(V[5]>R\),
take the inverse image in \(V[5]\) of an invariant line in
\(V[5]/R\). It is isotropic since \(\beta(R,V[5])=0\) and it is
generated by \(R\) and one further element. The only remaining
possibility is \(V\cong C_{25}\times C_{25}\). For generators
\(x,y\) of the two factors the entire pairing takes values in
\(\langle\beta(x,y)\rangle\), which is cyclic of order \(25\) by
nondegeneracy. The order formula for a perfect cyclic-valued
pairing then bounds the order of an isotropic subgroup by \(25\),
inconsistent with the given witness of order at least \(5^3\).
This proves Proposition~\ref{prop:relative} in the remaining case.
\end{proof}

\section*{Acknowledgements}
This work was supported by the National Natural Science Foundation of
China (grant 12371021).


\begin{thebibliography}{9}

\bibitem{AG}
J.~L. Alperin and G. Glauberman,
\href{https://math.uchicago.edu/~gg/limits.pdf}{Limits of abelian subgroups of finite \(p\)-groups},
\emph{J. Algebra} \textbf{203} (1998), 533--566.

\bibitem{Ber}
Y. Berkovich,
\href{https://doi.org/10.1006/jabr.2001.8790}{On subgroups and epimorphic images of finite \(p\)-groups},
\emph{J. Algebra} \textbf{248} (2002), 472--553.

\bibitem{Gl}
G. Glauberman,
\href{https://doi.org/10.1090/conm/524/10345}{A \(p\)-group with no normal large abelian subgroup},
in \emph{Character theory of finite groups}, Contemp. Math. \textbf{524},
Amer. Math. Soc., Providence, RI, 2010, 61--65.

\bibitem{GlExist}
G. Glauberman,
\href{https://doi.org/10.1016/j.jalgebra.2006.08.037}{Existence of normal subgroups in finite \(p\)-groups},
\emph{J. Algebra} \textbf{319} (2008), 800--805.

\bibitem{KN}
E.~I. Khukhro and V.~D. Mazurov (eds.),
\href{https://alglog.org/21tkt.pdf}{\emph{The Kourovka Notebook. Unsolved problems in group theory}},
21st ed., 2026.

\bibitem{KJ}
M. Konvisser and D. Jonah,
\href{https://doi.org/10.1016/0021-8693(75)90186-6}{Counting abelian subgroups of \(p\)-groups. A projective approach},
\emph{J. Algebra} \textbf{34} (1975), 309--330.

\end{thebibliography}
\end{document}